\documentclass[runningheads]{llncs}
\usepackage[T1]{fontenc}
\usepackage{graphicx}
\usepackage[noadjust,compress]{cite}

\usepackage{algorithm}
\usepackage{algpseudocode}
\usepackage{amsmath}
\usepackage{enumerate}
\usepackage{multirow}
\usepackage[normalem]{ulem}

\usepackage{hyperref}

\DeclareMathOperator{\NN}{NN}

\usepackage{amsmath,amsfonts,bm}

\def\eqref#1{equation~\ref{#1}}

\def\1{\bm{1}}

\def\vb{{\bm{b}}}
\def\vc{{\bm{c}}}
\def\vd{{\bm{d}}}

\def\vg{{\bm{g}}}
\def\vh{{\bm{h}}}

\def\vx{{\bm{x}}}
\def\vy{{\bm{y}}}
\def\vz{{\bm{z}}}

\def\evg{{g}}
\def\evh{{h}}

\def\evx{{x}}
\def\evy{{y}}

\def\mW{{\bm{W}}}

\DeclareMathAlphabet{\mathsfit}{\encodingdefault}{\sfdefault}{m}{sl}
\SetMathAlphabet{\mathsfit}{bold}{\encodingdefault}{\sfdefault}{bx}{n}

\def\sL{{\mathbb{L}}}

\def\sN{{\mathbb{N}}}

\def\sS{{\mathbb{S}}}

\def\sX{{\mathbb{X}}}

\begin{document}
\title{Optimization over Trained (and Sparse) Neural Networks: A Surrogate within a Surrogate}
\titlerunning{Optimization over Trained (and Sparse) Neural Networks}
%
\author{Hung Pham\inst{1} \and
Aiden Ren\inst{1} \and
Ibrahim Tahir\inst{1} \and
Jiatai Tong\inst{2} \and
Thiago Serra\inst{3}}
\authorrunning{H. Pham et al.}
%
\institute{Bucknell University, Lewisburg PA, United States\\
\email{\{hqp001,zr002,it005\}@bucknell.edu}
\and 
Northwestern University, Evanston IL, United States\\
\email{jiataitong2026@u.northwestern.edu}
 \and
University of Iowa, Iowa City IA, United States\\
\email{thiago-serra@uiowa.edu}}
\maketitle              
\begin{abstract}
In constraint learning, 
we use a neural network as a \emph{surrogate} for 
part of the constraints or of the objective function of 
an optimization model.
However, 
the tractability of the resulting model is heavily influenced by 
the size of the neural network used as a surrogate. 
One way to obtain a more tractable surrogate is by pruning the neural network first.
In this work, 
we consider how to approach the setting in which 
the neural network is actually a given: 
how can we solve an optimization model embedding a large and predetermined neural network?
We propose surrogating the neural network itself by pruning it, 
which leads to a sparse and more tractable optimization model, 
for which we hope to still obtain good solutions with respect to the original neural network. 
For network verification and function maximization models, 
that indeed leads to better solutions within a time limit, 
especially---and surprisingly---if we skip the standard retraining step known as finetuning. 
Hence, 
a pruned network with worse inference for lack of finetuning 
can be a better surrogate.


\end{abstract}

\section{Introduction}\label{sec:introduction}

In the last five years, we have seen a growing interest in 
approximating a constraint or an objective function of an optimization model with a neural network.
This approach is often denoted as \emph{constraint learning}. 
The most compelling circumstance for using constraint learning is when the exact form of some constraints or part of the objective function 
is unknown, but can be approximated using available data. 
When the exact form is known, 
another compelling circumstance is when the formulation is intractable for a combination of factors such as being nonlinear, nonconvex, and very large~\cite{misener2023surrogate,wang2023ensemble,chi2024trust,zhu2024validity}. 

Constraint learning has been applied in optimization models for  
scholarship allocation~\cite{bergman2022janos},  chemotherapy~\cite{maragno2023mixed}, molecular design~\cite{mcdonald2023molecular}, power grid operation~\cite{chen2020voltage,murzakhanov2022powerflow,liu2025doc}, 
and automated control in general~\cite{say2017planning,wu2020scalable,yang2021control,sosnin2024control}. 
We can also use optimization models involving neural networks to evaluate those neural networks  
for adversarial perturbations~\cite{cheng2017resilience,anderson2020strong,rossig2021verification,strong2021global,depalma2021convex}, 
compression~\cite{serra2020lossless,serra2021compression,elaraby2023oamip}, counterfactual explanations~\cite{kanamori2021counterfactual,tsiourvas2024counterfactual}, equivalence~\cite{kumar2019equivalent}, expressiveness~\cite{serra2018bounding,serra2020empirical,cai2023pruning}, 
monotonicity~\cite{liu2020monotonic}, 
and reachability~\cite{LomuscioMIP,DuttaMIP}. 
Conversely, 
when the neural network is trained as a reinforcement learning policy, 
the constraints in the optimization model 
can be used for limiting the action space~\cite{delarue2020rlvrp,burtea2023safe}. 
Many frameworks have been proposed 
to embed neural networks and other machine learning models as part of optimization models, 
such as JANOS~\cite{bergman2022janos}, reluMIP~\cite{lueg2021relumip}, OMLT~\cite{ceccon2022omlt}, OCL~\cite{fajemisin2023ocl}, OptiCL~\cite{maragno2023mixed}, Gurobi Machine Learning~\cite{gurobi2025ml}, and PySCIPOpt-ML~\cite{turner2023pyscipoptml}.

Some neural network architectures are more convenient to embed in optimization models than others. 
In particular, 
we typically use the Rectified Linear Unit (ReLU) activation function~\cite{hahnloser2000origin,nair2010rectified} for a couple of reasons. 
First, 
this activation function became widely popular in the 2010s due to its good classification performance and relatively lower training  cost~\cite{nair2010rectified,glorot2011rectifier,lecun2015nature}. 
Similar to what was already known for other neural network architectures much earlier~\cite{cybenko1989approximation,funahashi1989approximate,hornik1989approximator}, 
we now know that ReLU networks are also universal function approximators~\cite{yarotsky2017relu,hanin2017approximating}.  
More recently, Gaussian Error Linear Units (GELUs)~\cite{hendrycks2016gelus} have been used from the very beginning in many foundation models based on the transformer architecture~\cite{vaswani2017attention}, 
such as GPT~\cite{radford2018gpt}, BERT~\cite{devlin2019bert}, and ViT~\cite{dosovitskiy2021vit}. 
Nevertheless, a ReLU is a piecewise linear approximation of a GELU, 
which brings us to the next point. 
Second, 
and perhaps most importantly for the continued use of ReLUs nowadays, 
each neuron with ReLU activation represents a piecewise linear function---and by consequence a neural network with only ReLUs is also a piecewise linear function~\cite{arora2018understanding}. 
In the context of nonlinear optimization, 
substantial work has been devoted to using piecewise linear approximations in surrogate models~\cite{mccormick1976under,chien1977simplicial,rosen1986separable,feijoo1988nonsep,babayeve1997two,mangasarian2005under,magnani2009convex,vielma2010framework,misener2010multi,williams2013model,huchette2023pwl} 
because we can resort to linear optimization for iterative improvements over such approximations. 
In fact, there are many active lines of research on what piecewise linear representations can be obtained from different neural network architectures, most of which focused on ReLUs~\cite{pascanu2013on,montufar2014on,telgarsky2015benefits,montufar2017notes,raghu2017expressive,arora2018understanding,serra2018bounding,hinz2019framework,hanin2019complexity,hanin2019deep,serra2020empirical,phuong2020equivalence,rolnick2020reverse,hertrich2021depth,tseran2021expected,wang2022estimation,grigsby2022topology,montufar2022maxout,cai2023pruning,haase2023lower,huchette2023survey,safran2024max,hertrich2024virtual,averkov2025depth}.

With a piecewise linear representation, 
we can embed the neural network in the optimization model using a Mixed-Integer Linear Programming~(MILP) formulation. 
However, such formulations range from being small but having a weak linear relaxation to having a stronger relaxation but being prohibitivelly large~\cite{huchette2018phd}, 
which makes it difficult to use an off-the-self MILP solver. 
That motivated studies on how to tackle such models, ranging from  reformulation~\cite{huchette2018phd,schweidtmann2019global,tsay2021partition,liu2025doc} to methods for obtaining smaller big M coefficients through activation bounds~\cite{cheng2017resilience,fischetti2018constraints,grimstad2019surrogate,liu2021algorithms,tsay2021partition,badilla2023tradeoff,zhao2024horizon,hojny2024message,sosnin2024control} and parameter rescaling~\cite{plate2025scaling},
activation inferences for fixing variables and limiting search space~\cite{tjeng2019evaluating,xiao2019training,botoeva2020efficient,serra2020empirical}, cutting planes~\cite{anderson2020strong},  and heuristic solutions~\cite{perakis2022optimizing,tong2024walk,tong2026gradient}. 
Another---perhaps overlooked---approach is using neural networks that are smaller~\cite{cacciola2024structured} and sparser~\cite{say2017planning,xiao2019training}.

In this work, 
we explore this latter approach of using sparser neural networks. 
On the one hand, 
we know that (i) sparser optimization models tend to be solved faster; 
and (ii) we can sparsify a neural network, which is known as \emph{network pruning}, and still recover its performance by retraining for a few steps, which is known as \emph{finetuning}. 
In fact, 
it is already reasonably well established that pruning a neural network in advance leads to a more tractable constraint learning model~\cite{say2017planning,xiao2019training}. 
On the other hand, 
finetuning requires access to training data and requires a non-negligible runtime overhead. 
Moreover,  
if the purpose of solving an optimization model that embeds a neural network is to verify properties of the neural network itself, 
then it is not immediate if---and how---to use a pruned counterpart of the given neural network as a proxy. 
Therefore, 
we evaluate how much we can prune, and how much we should finetune, a neural network to still obtain effective and efficient surrogate optimization models.

\section{From Notation to Formulation}

In this paper, we consider feedforward networks with fully-connected layers of neurons having ReLU activation. 
Note that convolutional layers can be represented as fully-connected layers with a block-diagonal weight matrix, 
for which reason we abstract that possibility. 
We also abstract that fully-connected layers are often followed by a softmax layer~\cite{bridle1990softmax}, 
since the largest output of softmax matches the largest input of softmax, 
for which reason it is not necessary to include the softmax layer in applications such as network verification.

Each neural network has an input $\vx = [\evx_1 ~ \evx_2 ~ \dots ~ \evx_{n_0}]^\top$ from a bounded domain $\sX$ and corresponding output $\vy = [\evy_1 ~ \evy_2 ~ \dots ~ \evy_m]^\top$, and each layer $l \in \sL = \{1,2,\dots,L\}$ has output $\vh^l = [\evh_1^l ~ \evh_2^l \dots \evh_{n_l}^l]^\top$ from neurons indexed by $i \in \sN_l = \{1, 2, \ldots, n_l\}$. 
Let $\mW^l$ be the $n_l \times n_{l-1}$ matrix where each row corresponds to the weights of a neuron of layer $l$, $\mW_i^l$ the $i$-th row of $\mW^l$, and $\vb^l$ the vector of biases associated with the units in layer $l$. With $\vh^0$ for $\vx$ and $\vh^{L}$ for $\vy$,  the output of each unit $i$ in layer $l$ consists of an affine function $\evg_i^l = \mW_{i}^l \vh^{l-1} + \vb_i^l$ followed by the ReLU activation $\evh_i^l = \max\{0, \evg_i^l\}$. 
When training the neural network, 
we vary the paramenters in $\{ (\mW^l, \vb^l) \}_{l \in \sL}$ 
to better fit the values given for the set of input--output pairs $\{ (\vx^{(i)}, \vy^{(i)}) \}_{i \in \sS}$ representing the training set.

When optimizing over the trained neural network, 
we flip what is variable and what is constant in relation to training: 
we vary the input $\vx = \vh^0$ and the outputs at each layer $\{ (\vg^l, \vh^l) \}_{l \in \sL}$ for the fixed parameters $\{ (\mW^l, \vb^l) \}_{l \in \sL}$. 
For each layer $l \in \sL$ and neuron $i \in \sN_l$, 
we also use a binary variable $\vz^l_i$ representing the ReLU activation to map inputs to outputs using MILP:
\begin{align}
    \mW_i^l \vh^{l-1} + \vb_i^l = \vg_i^l \label{eq:mip_unit_begin} \\
    (\vz_i^l = 1) \rightarrow (\vh_i^l = \vg_i^l) \label{eq:first_indicator} \\  
    (\vz_i^l = 0) \rightarrow (\vg_i^l \leq 0 \wedge \vh_i^l = 0) \label{eq:last_indicator} \\
    \vh_i^l \geq 0 \label{eq:lp_unit_end} \\ 
    \vz^l_i \in \{0,1\} \label{eq:mip_unit_end}
\end{align}
The indicator constraints (\ref{eq:first_indicator})--(\ref{eq:last_indicator}) can be modeled with big M constraints~\cite{bonami2015indicator}.

In the absence of other decision variables, 
one general form of representing a linear optimization model embedding a neural network is as follows:
\begin{subequations}
\begin{align}
\max ~~& \vc^T \vx + \vd^T \vy \\
\text{s.t.} ~~& A \vx + B \vy \leq b \\
& \vy = \NN(\vx)
\end{align}
\end{subequations}
We use $\vy = \NN(\vx)$ as a shorthand for the input--output mapping defined by the set of constraints $(\ref{eq:mip_unit_begin})\text{--}(\ref{eq:mip_unit_end}) ~ \forall l \in \sL, i \in \sN_l$ across the entire neural network.

Next, we describe two applications based on special cases of this formulation. 
They are both used for evaluating our approach.

\subsection{Network Verification}\label{sec:nv}

For neural networks used for classification, 
one application of embedding them in optimization models is to determine if there is an adversarial perturbation for a chosen input 
$\vx^{(i)}, i \in \sS$ from the training set. 
If $\vx^{(i)}$ is correctly classified as class $j \in \{1, \ldots, m\}$ by having an output $\vy$ such that $\evy^{(i)}_j > \evy^{(i)}_k ~ \forall k \in \{1, \ldots, m\} \setminus \{ j \}$, 
then we can try to determine if there is a similar input $\vx$ 
with a different classification. 
By varying the inputs within $\{ \vx : \| \vx - \vx^{(i)} \|_{1} \leq \varepsilon \}$ for a chosen $\varepsilon$,  
we try to find an input that is better classified with a chosen class $j' \neq j$, 
i.e., $\evy_{j'} > \evy_j$.\footnote{With $j'$ fixed, we are only ensuring that $j'$ would be a better classification than $j$, since there might be another class $j''$ dominating both, i.e.,  $\evy_{j''} > \evy_{j'} > \evy_j$.}
That leads to the following MILP formulation: 
\begin{subequations}
\begin{align}
\max ~~& \evy_{j'} - \evy_j \label{eq:vnn_begin} \\
\text{s.t.} ~~& \sum_{k \in \{1, \ldots, n_0\}} |\vx_k - \vx^{(i)}_k| \leq \varepsilon   
\\& \vy = \NN(\vx)  \label{eq:vnn_end}
\end{align}
\end{subequations}
In the model above, 
any solution with a positive objective function value entails an adversarial perturbation; 
whereas 
a nonpositive optimal value implies that no such perturbation exists 
for the given choices of $i$, $j'$, and $\varepsilon$.

\subsection{Function Maximization}\label{sec:fm}

For neural networks used for regression, 
one application of embedding them in optimization models is to 
optimize over the function surrogated by the neural network. 
The case of finding the maximum for a neural network with a single output---i.e., $m=1$---leads to the following MILP formulation:
\begin{subequations}
\begin{align}
\max ~~& \evy_1 \label{eq:fm_begin} \\
\text{s.t.} ~~& \vy = \NN(\vx) \\
& \vx \in \sX \label{eq:fm_end}
\end{align}
\end{subequations}

\section{From Network Pruning to Sparse Surrogates}

Neural networks have a very peculiar trait: 
assuming that two neural networks can be trained to achieve the same level of accuracy, 
it is often easier to train the largest one than it is to train the smallest one. 
But after training a neural network that is larger than it needs to be, 
we can simplify the network by removing neurons or connections and then still recover a similar accuracy by carefully adjusting the remaining parameters. 
In fact, 
this is becoming mainstream knowledge
with the constant discussion about number of parameters in large language models and the application of pruning techniques for obtaining comparable variants that are smaller and faster for inference~\cite{frantar2023sparsegpt,ma2023pruner,sun2024pruning,xia2024sheared}.

\subsection{Background}

Neural networks are pruned by either (i) removing connections, 
which is equivalent to zeroing out specific parameters---known as \emph{unstructured} pruning; 
or (ii) removing units, such as neurons, convolutional filters, or layers---known as \emph{structured} pruning. 
The latter has greater appeal for performance, 
since it goes beyond reducing  storage 
to using smaller hardware and running the model faster, 
but then we need to prune less to still recover the original performance~\cite{cheng2018compression}.

\textbf{But why do we need network pruning?}  
A larger neural network has a smoother loss landscape~\cite{li2018landscape,ruoyu2020landscape}, 
which facilitates training convergence; 
and the larger size may also prevent layers from becoming inactive~\cite{riera2022jumpstart}, 
which is a common cause for unsuccessful training. 
\textbf{Why does network pruning work?} 
In larger networks,
there is redundancy among the parameters~\cite{denil2013parameters}, 
and zeroing out parameters leads to a loss landscape from which finetuning the pruned network for recovering the original performance converges considerably faster than the original training~\cite{jin2022generalization}. 
From a model flexibility perspective, 
unstructured pruning at moderate rates has little effect on the expressiveness of the neural network architecture~\cite{cai2023pruning}.
\textbf{How much can we prune?} 
In sufficiently large networks, 
we can remove as much as half of the parameters and still recover the original performance---or even improve upon it~\cite{hoefler2021sparsity}. 
However, that varies with the task for which the network is  trained~\cite{liebenwein2021lost}. 
Moreover, 
pruning may have a disparate effect across classes~\cite{hooker2019forget,paganini2020responsibly,hooker2020bias}, 
which may lead to pruned networks that exacerbate existing performance differences~\cite{tran2022disparate}. 
On the bright side, 
a smaller amount of pruning may actually correct such distortions~\cite{good2022recall}. 
\textbf{What should we prune?} 
The two main philosophies~\cite{blalock2020survey} are (1) to remove parameters with the smallest absolute value---dating back to \cite{hanson1988minimal,mozer1989relevance,janowsky1989prunning}; and (2) to remove parameters with the smallest expected impact on the output---dating back to \cite{lecun1989damage,hassibi1992surgeon,hassibi1993surgeon}, and including the special case of exact compression~\cite{serra2020lossless,sourek2021lossless,serra2021compression,ganev2022compression}.  
\textbf{And when should we prune?} 
Most studies have focused on pruning once (\emph{one-shot}) and after training, 
but 
recent work has shown that it might be beneficial to prune iteratively and during training~\cite{frankle2019lottery}, 
or even before training~\cite{lee2019pretraining}.

Mathematical optimization has been extensively used in more sophisticated pruning methods~\cite{lecun1989damage,hassibi1992surgeon,hassibi1993surgeon,he2017cnn,luo2017thinet,aghasi2017nettrim,serra2020lossless,ye2020forward,singh2020woodfisher,verma2021subdifferential,serra2021compression,yu2022cbs,cai2023pruning,elaraby2023oamip,benbaki2023chita,ebrahimi2023kfac,cacciola2023perspective,wu2024iobs}. 
Those methods tend to perform better than simpler heuristics at higher pruning rates. 
However, they also come at a greater computational cost. 
For moderate pruning, something as simple as removing the weights with the smallest absolute values, known as \emph{Magnitude Pruning} (MP)~\cite{hanson1988minimal,mozer1989relevance,janowsky1989prunning}, 
remains a very competitive choice~\cite{yu2022cbs}. 
We take the liberty of calling it ``unreasonably effective''.

The use of pruned neural networks in mathematical optimization, 
however, 
has been less explored. 
Among the first studies on embedding neural networks, 
Say et al.~\cite{say2017planning} observed that a modest amount of unstructured pruning---removing about 20\% of the parameters---significantly reduced the runtime for solving the optimization model. 
Moreover, 
Xiao et al.~\cite{xiao2019training} and 
Cacciola et al.~\cite{cacciola2024structured} 
observed that structured pruning---having fewer neurons and therefore fewer binary decision variables mapping the activation state of each neuron---leads to comparable neural networks that are more easily verified.

\subsection{The Sparse Surrogate Approach}

Suppose that we have a (dense) neural network $\mathcal{D}$ to which $\mathcal{S}$ is a sparse counterpart obtained by network pruning, 
with $\vy^{D} = \mathcal{D}(\vx)$ and $\vy^{S} = \mathcal{S}(\vx)$ as the corresponding outputs from those two neural networks for a same input $\vx$.

We will succinctly describe how to use $\mathcal{S}$ for obtaining solutions for constraint learning models on $\mathcal{S}$. We will use the models from Sections~\ref{sec:nv} and \ref{sec:fm}.

\subsubsection{Network Verification}

First, we consider the case of a network verification problem, in which we validate if an adversarial input to the pruned network is an adversarial input to the original network. 
Let $\text{VNN}(\mathcal{N}, \vx^{(i)}, \varepsilon, j, j')$ be the MILP formulation (\ref{eq:vnn_begin})--(\ref{eq:vnn_end}) for verifying neural network $\mathcal{N} \in \{ \mathcal{D}, \mathcal{S} \}$ starting from input $\vx^{(i)}$ and with maximum $L1$-norm distance $\varepsilon$ for obtaining another input $\vx$ in which the output for class $j'$ is as large as possible in comparison to that of class $j$, i.e.,
we want to find an input $\vx$ maximizing $\evy^N_{j'} - \evy^N_{j}$ for $\vy^N = \mathcal{N}(\vx)$. \linebreak
For the purpose of verification, 
it suffices to find an $\vx$ such that $\evy^N_{j'} > \evy^N_{j}$. 
To find a solution with positive value for dense model $\textbf{VNN}(\mathcal{D}, \varepsilon, x^{(i)}, j, j')$, 
we resort to solving sparse model $\textbf{VNN}(\mathcal{S}, \varepsilon, x^{(i)}, j, j')$ as outlined in Algorithm~\ref{alg:callback}.

\begin{algorithm}[b]
    \caption{Heuristic for obtaining an adversarial input to a (dense) neural network $\mathcal{D}$ by trying to solve the same problem on its sparse counterpart $\mathcal{S}$} 
    \label{alg:callback}
    {
    \begin{algorithmic}[1]
    \While{trying to solve $\textbf{VNN}(\mathcal{S}, \varepsilon, x^{(i)}, j, j')$} \label{lin:while} \Comment{MILP solver call}
        \If{solution $(\vx, \vy^S)$ found} \Comment{Feasible solution callback} \label{lin:callback}
            \State $\vy^D = \mathcal{D}(\vx)$ \Comment{Output of the dense model $\mathcal{D}$} \label{lin:yD}
            \If{$\evy_{j'}^D > \evy_j^D$} \Comment{Check if $\vx$ is adversarial to $\mathcal{D}$} \label{lin:check}
                \State \Return $\vx$ \Comment{Adversarial input found}
            \EndIf
        \EndIf
    \EndWhile
    \State \Return $\emptyset$ \Comment{No adversarial input found}
    \end{algorithmic}
    }
\end{algorithm}


\subsubsection{Function Maximization} Next, we consider the case of a function maximization problem, in which we check if an input to the pruned network yields a better value than the inputs previously tried on the dense network. 
Let $\textbf{FM}(\mathcal{N})$ be the MILP formulation (\ref{eq:fm_begin})--(\ref{eq:fm_end}) for maximizing the output of the function modeled by network $\mathcal{N} \in \{ \mathcal{D}, \mathcal{S} \}$ over domain $\sX$. 
Unlike the network verification case, there is no criterion for prematurely stopping the search in this application.
To find a solution with better value for dense model $\textbf{FM}(\mathcal{D})$, 
we resort to solving sparse model $\textbf{FM}(\mathcal{S})$ as outlined in Algorithm~\ref{alg:fm_callback}. 

\begin{algorithm}[h]
    \caption{Heuristic for maximizing the output of a (dense) neural network $\mathcal{D}$ by trying to solve the same problem on its sparse counterpart $\mathcal{S}$} 
    \label{alg:fm_callback}
    {
    \begin{algorithmic}[1]
    \State $\vx^* \gets \emptyset$ \Comment{Initialize best solution as none}
    \State $y^* \gets - \infty$ \Comment{Placeholder for best solution value}
    \While{trying to solve $\textbf{FM}(\mathcal{S})$} \label{lin:fm_while} \Comment{MILP solver call}
        \If{solution $(\vx, y^S)$ found} \Comment{Feasible solution callback} \label{lin:fm_callback}
            \State $y^D = \mathcal{D}(\vx)$ \Comment{Output of the dense model $\mathcal{D}$} \label{lin:fm_yD}
            \If{$\vx^* = \emptyset$ or $y^D > y^*$} \Comment{Check if $\vx$ is the first or a better solution} \label{lin:fm_check}
                \State $\vx^* \gets \vx$ \Comment{Update best solution}
                \State $y^* \gets y^D$ \Comment{Update best solution value}
            \EndIf
        \EndIf
    \EndWhile
    \State \Return $\vx^*$ \Comment{Provide best solution found}
    \end{algorithmic}
    }
\end{algorithm}


In both cases, we assume that the algorithms will be run up to a reasonable time limit, 
so that we would not be done solving the sparse model but would have obtained and evaluated a sufficient number of solutions in the dense model.

\section{Experiments for Network Verification}\label{sec:nv_exp}

We evaluated the time for finding an adversarial input for a (dense) neural network $\mathcal{D}$ by directly solving model $\text{VNN}(\mathcal{D}, \vx^{(i)}, \varepsilon, j, j')$, which we denote as \emph{Dense Runtime},  
in comparison to indirectly solving model $\text{VNN}(\mathcal{S}, \vx^{(i)}, \varepsilon, j, j')$ 
while resorting to Algorithm~\ref{alg:callback}, 
which we denote as \emph{Pruned Runtime}.
Our goal is to find if, and when, Pruned Runtime is shorter than Dense Runtime.

\subsection{Technical Details}

We used the source code of SurrogateLIB~\cite{turner_2024_11231147} as the basis for training neural networks and producing network verification problems, 
starting with the MNIST dataset~\cite{lecun1998mnist} and then extending the code to also work with the Fashion-MNIST dataset~\cite{xiao2017fashion}.
For each of those datasets, 
we tried all combined variations of inputs sizes $n_0 = 18^2$ (compressed images) and $n_0 = 28^2$ (images at original size), 
number of ReLU layers $L \in \{2, 4\}$, 
and uniform layer width $n_i \in \{32, 64\} ~\forall i \in \{1, \ldots, L\}$. 
For each dataset and choice of hyperparameters, 
we used 10 randomization seeds for training neural networks,  associating them with verification problems on distinct samples from the training set, 
and choosing some $\varepsilon \in [4.5, 5.5]$. 
In total, we have 160 verification problems.

For producing pruned versions of those networks, 
we used the PyTorch library. 
On the number of parameters removed, we applied layerwise pruning rates of 0.3, 0.5, 0.8, 0.9, and 0.95.
On what to prune, 
we applied both Magnitude Pruning~(MP), 
which corresponds to pruning the parameters with the smallest absolute value; 
and Random Pruning~(RP), 
which corresponds to pruning parameters randomly. 
In addition, we also evaluated \emph{unstructured pruning}, in which case the parameters across the layer are pruned indiscriminately; and \emph{structured pruning}, in which case we consider all the parameters associated with a neuron and decide about pruning whole neurons instead. 
Finally, 
as a last step we also opted between finetuning the pruned network or not finetuning it and keeping it as it was after being pruned.
For each of the 160 original verification problems, 
the combination of all the pruning choices above resulted in solving variants in 40 pruned versions of each neural network. 
In total, we use 6,560 models.

We used the BisonNet cluster. The steps involving the solution of MILP models with Gurobi were run on Intel Xeon Gold 6442Y CPUs. The steps involving network training and pruning were run on AMD EPYC 7252 CPUs with NVIDIA RTX A5000 GPUs. 
The neural networks were trained and pruned using Torch 2.0.0. 
Each model on MNIST was trained for 5 epochs, 
and each model on Fashion-MNIST was trained for 40 epochs.
Without finetuning, 
there was a single round of pruning. 
With finetuning,
there were 5 pruning rounds and each round had 5 epochs of retraining. 
The MILP models were solved using Gurobi Optimizer 10.0.1 
with a time limit of 300 seconds for each of the 6,560 models.

\subsection{Results and Analysis}

\begin{figure}[h!]
\centering
\includegraphics[height=0.85\textheight]{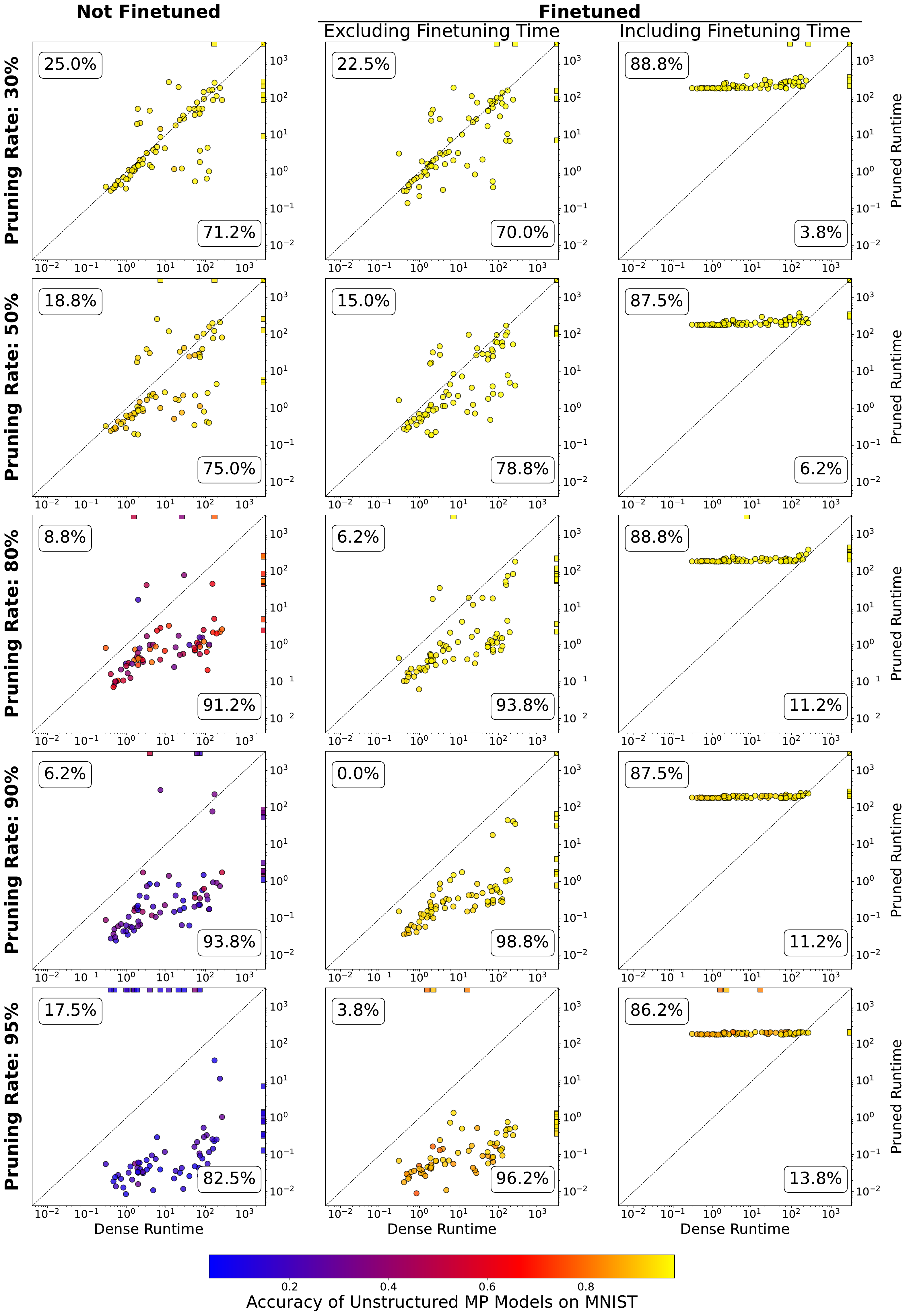}
\caption{Time to find adversarial input to networks trained on MNIST by solving the verification problem directly ($x$ axis) or indirectly with Algorithm~\ref{alg:callback} ($y$ axis) per pruning rate, use of finetuning, and inclusion of finetuning in runtime.
Squares on top or (and) right sides indicate no adversarial input found for either (both). 
Ties are not counted.
}\label{fig:prune_mnist}
\end{figure}

\begin{figure}[h!]
\centering
\includegraphics[height=0.85\textheight]{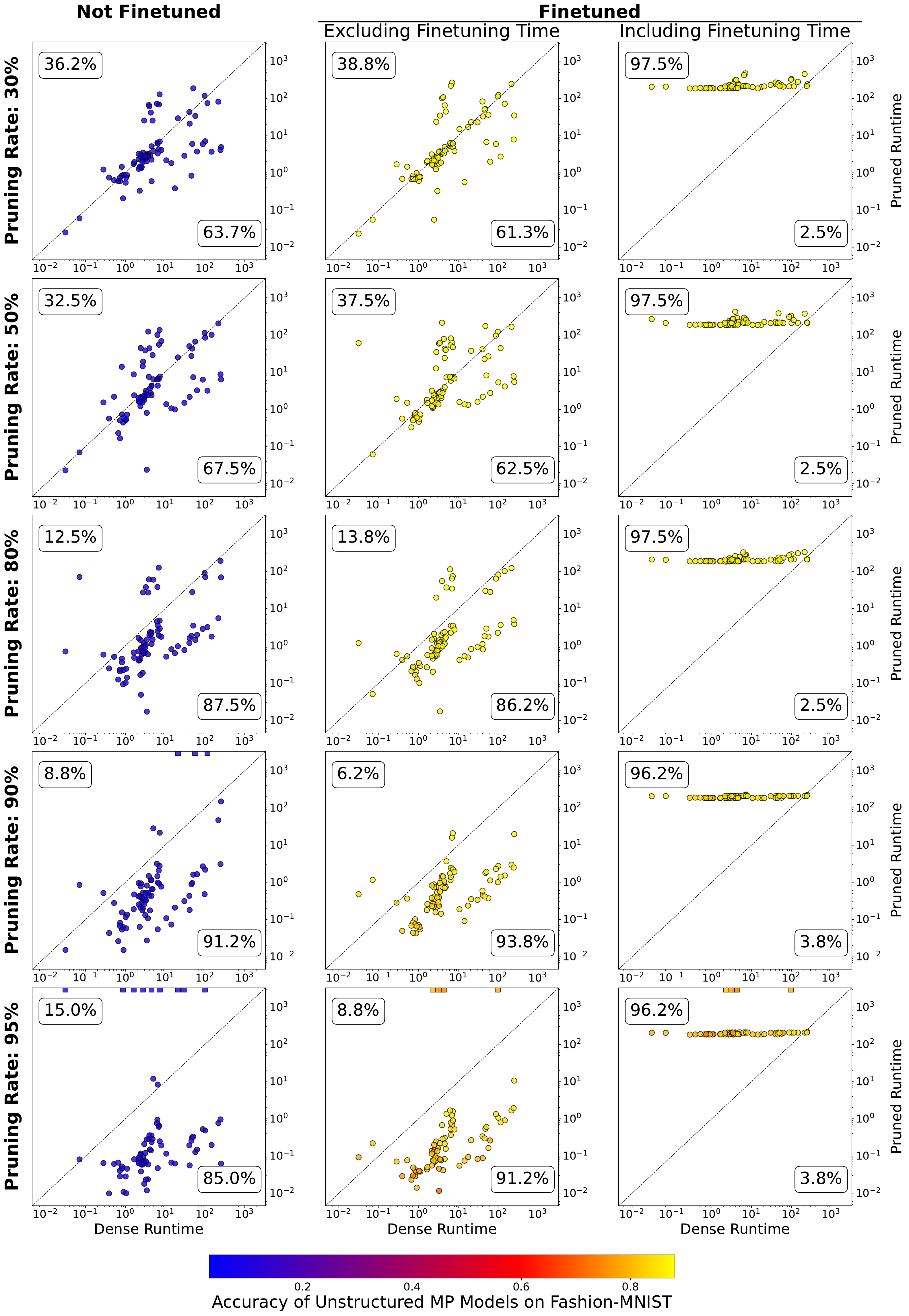}
\caption{Time to find adversarial input to networks on Fashion-MNIST by solving the verification problem directly ($x$ axis) or indirectly with Algorithm~\ref{alg:callback} ($y$ axis) per pruning rate, use of finetuning, and inclusion of finetuning in runtime.
Squares on top or (and) right sides indicate no adversarial input found for either (both). 
Ties are not counted.}\label{fig:prune_fashion}
\end{figure}

Our best results were obtained by using Algorithm~\ref{alg:callback} with unstructured MP instead of solving the verification model directly. 
We compare the runtimes (in seconds) to solve the verification model directly and indirectly for MNIST in Figure~\ref{fig:prune_mnist} and Fashion-MNIST in Figure~\ref{fig:prune_fashion}. 
We vary pruning rate, whether finetuning is used, and whether finetuning time is counted as part of the runtime.
We report the percentage of instances above and below the identity line, corresponding to the direct and the indirect approaches being faster. 
Those percentages do not add up to 100\% if there are ties, which includes when both methods time out.

We draw the following conclusions and provide a rationale to each conclusion by listing observations about the plots in those figures: 
\begin{enumerate}[(I)]
\item Using our approach in network verification is advantageous in terms of (i) individual runtimes as well as (ii) number of instances solved: 
\begin{enumerate}[(i)]
    \item We found adversarial inputs faster for most instances regardless of pruning rate and of finetuning the pruned neural networks or not. 
    \item When there are timeouts (i.e., not finding an adversarial input) from solving the verification problem directly (as with MNIST), then our approach with a small pruning rate reduces the number of timeouts.
\end{enumerate}

\item The pruning rate can be adjusted for different purposes, from finding more adversarial inputs up to a time limit at lower rates (ii again) to finding (iii) most adversarial inputs faster at higher rates and (iv) fewer adversarial inputs but in a much shorter amount of time at the highest rate: 
\begin{enumerate}[(i)]
    \setcounter{enumii}{2}
    \item The number of runtime improvements increases with pruning rate up to 90\% for both datasets, finetuned or not, and then decreases.
    \item The differences between runtimes become more extreme in our favor with greater pruning rates, but then the timeouts also increase. 
\end{enumerate}

\item The pruned neural network does not need to be a good classifier to help us find adversarial inputs (v). In fact, it is not helpful to finetune the network to improve accuracy after pruning at lower pruning rates (vi). 
For higher pruning rates, finetuning the neural network can be helpful (vii), 
but the cost of finetuning would have to be amortized over solving multiple verification problems on the same neural network (viii): 
\begin{enumerate}[(i)]
    \setcounter{enumii}{4}
    \item The lower accuracy in the case without finetuning, almost approaching random guessing (10\% on either dataset), did not prevent us from using the pruned neural networks for obtaining adversarial inputs. 
    \item The results were better without finetuning for the lowest pruning rates (one for MNIST and three for Fashion-MNIST).
    \item The percentage difference of instances solved faster with finetuning is only on the second most significant digit for the higher pruning rate, except the highest (e.g., 98.8\% instead of 93.5\% on MNIST and 93.8\% instead of 91.2\% on Fashion-MNIST for pruning rate 90\%).
    \item If we account for the finetuning cost, then it is generally faster to solve the verification problem directly. 
\end{enumerate}

\end{enumerate}

Given the cost--benefit advantage of not using finetuning, 
we conducted the following ablations restricted to the results without finetuning.

\setlength{\tabcolsep}{9pt}
\begin{table}[b]
\centering
\caption{Percentage of instances for which solving the verification problem indirectly on a pruned  network is faster by   
pruning rate (columns); 
using unstructured and structured pruning (top and bottom rows); and 
not finetuning (NF) or finetuning (F). 
Unfavorable figures (below 50\%) are reported in italics for greater emphasis.
}\label{tab:structured}
\begin{tabular}{l c || c c c c c}
  && \multicolumn{5}{c}{\textbf{Pruning Rate}} \\
  && 0.3 & 0.5 & 0.8 & 0.9 & 0.95 \\ 
  \hline \hline
  \multirow{2}{*}{\textbf{Unstructured}} & NF & 67.5\% & 71.3\% & 89.4\% & 92.5\% & 83.8\% \\
   & F & 65.6\% & 70.6\% & 90.0\% & 96.3\% & 93.8\% \\
   \hline
   \multirow{2}{*}{\textbf{Structured}} & NF & 75.0\% & 78.8\% & 70.6\% & 72.5\% & 65.0\% \\
   & F & 59.4\% & 57.5\% & 56.3\% & \emph{49.4\%} & \emph{46.3\%} \\
\end{tabular}
\end{table}

\begin{figure}
\centering
\includegraphics[height=0.85\textheight]{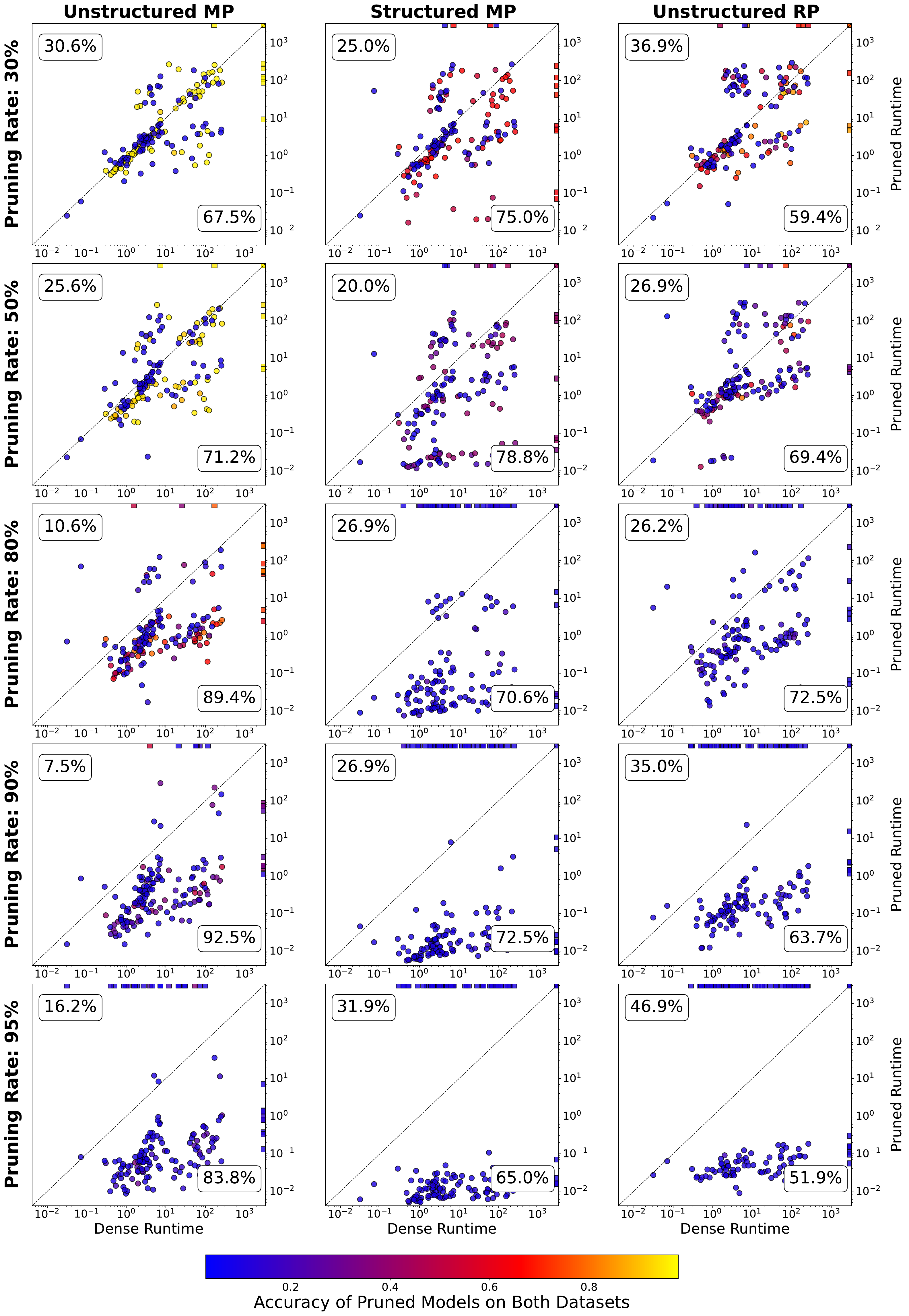}
\caption{Time to find adversarial input to networks on MNIST and Fashion-MNIST by solving the verification problem directly ($x$ axis) or indirectly with Algorithm~\ref{alg:callback} ($y$ axis) per pruning rate for different forms of network pruning. 
Squares on top or (and) right sides indicate no adversarial input found for either (both). 
Ties are not counted.
}\label{fig:ablations}
\end{figure}

First, we considered the impact of other network pruning choices. 
Figure~\ref{fig:ablations} compares the joint results on MNIST and Fashion-MNIST by using unstructured MP as before, then by only replacing unstructured with structured pruning, and then by only replacing MP with RP. 
Considering possible questions due to the favorable results for structured pruning in some of the plots, 
Table~\ref{tab:structured} summarizes the number of instances on both datasets in which an adversarial input would have been found faster than by directly solving the verification problem according to the structural type of pruning and the use of finetuning.

The purpose of considering RP, which is not an appealing choice intuitively, 
was to evaluate if a deliberate choice on how to prune would influence the results. 
We can see that it does---and that 
the only deliberate choice that we considered was MP. 
We opted for MP due to its small computational cost and widely regarded effectiveness in practice. 
Hence, we leave is open to future work if it would be beneficial to replace MP with a more sophisticated network pruning algorithm. 
That improvement seems reasonable to expect, but 
we believe that it would be beyond the scope of this particular study.

Second, we considered how the results vary based on the dimensions of the neural networks considered. 
Table~\ref{tab:dimensions} summarizes the number of instances on both datasets in which an adversarial input would have been found faster than by directly solving the verification problem according to input size, number of layers (network depth), and neurons per layer (layer width).

\setlength{\tabcolsep}{9pt}
\begin{table}[t]
\centering
\caption{Percentage of instances for which solving the verification problem indirectly on a pruned  network is faster by   
pruning rate (columns); 
disaggregated in terms of input size (top rows), number of layers (middle), and neurons per layer (bottom rows).}\label{tab:dimensions}
\begin{tabular}{l c || c c c c c}
  && \multicolumn{5}{c}{\textbf{Pruning Rate}} \\
  && 0.3 & 0.5 & 0.8 & 0.9 & 0.95 \\ 
  \hline \hline
  \multirow{2}{*}{\textbf{Input Size}} & $18^2$ & 75.0\% & 80.6\% & 86.9\% & 84.4\% & 76.3\% \\
   & $28^2$ & 67.5\% & 69.4\% & 73.1\% & 80.6\% & 72.5\% \\
   \hline
   \multirow{2}{*}{\textbf{Network Depth}} & 2 & 80.6\% & 84.4\% & 86.3\% & 86.3\% & 76.3\% \\
   & 4 & 61.9\% & 65.6\% & 73.8\% & 78.8\% & 72.5\% \\
   \hline
   \multirow{2}{*}{\textbf{Layer Width}} & 32 & 69.4\% & 78.8\% & 78.1\% & 81.9\% & 75.6\% \\
   & 64 & 73.1\% & 71.3\% & 81.9\% & 83.1\% & 73.1\% \\
\end{tabular}
\end{table}

We draw more conclusions and justify them based on the above ablations:
\begin{enumerate}[(I)]
\setcounter{enumi}{3}
\item Structured pruning without finetuning is also potentially applicable at lower pruning rates (ix), 
with the caveat that more timeouts may occur (x). 
However, 
finetuning after structured pruning 
appears to make the pruned networks significantly different from the original neural network, 
since their adversarial inputs are less compatible (xi):
\begin{enumerate}[(i)]
    \setcounter{enumii}{8}
    \item For the two lowest pruning rates (up to 50\%), more instances are solved faster with structured MP than with unstructured MP.
    \item Across all pruning rates, structured pruning leads to more timeouts.
    \item Finetuning has a positive effect when applied to neural networks that were subject to unstructured pruning, 
and that positive effect grows with pruning rate. The opposite happens with structured pruning.
\end{enumerate}

\item
Pruning criteria (such as MP vs. RP) affect the 
results significantly (xii), and 
may help tackling networks with larger inputs and more layers (xiii):
\begin{enumerate}[(i)]
    \setcounter{enumii}{11}
    \item The results for unstructured MP are better than those for unstructured RP, 
and the difference grows with the pruning rate.
    \item The benefit is more significant in neural networks with smaller input size or smaller depth, but less affected by layer width. 
\end{enumerate}

\end{enumerate}

\section{Experiments for Function Maximization}

We evaluated the best solution obtained for a (dense) neural network $\mathcal{D}$ by directly solving model $\textbf{FM}(\mathcal{D})$, 
which we denote as \emph{Dense Maximum Value}, 
in comparison to indirectly solving model $\textbf{FM}(\mathcal{S})$ while resorting to Algorithm~\ref{alg:fm_callback}, 
which we denote as \emph{Pruned Maximum Value}. 
Our goal is to find if, and when, Pruned Maximum Value is greater than Dense Maximum Value.

\subsection{Technical Details}

We used randomly initialized networks with 5 seeds for all combinations of 
input sizes $n_0 \in \{100, 1000, 10000\}$, number of layers $L \in \{2, 3, 4, 5\}$, uniform layer width $n_i \in \{50, 100, 200\} ~\forall i \in \{1, \ldots, L\}$, and with an extra one-neuron output layer. 
That setting is similar to other papers using this problem~\cite{perakis2022optimizing,tong2024walk}. 
Based on the results from Section~\ref{sec:nv_exp}, we used only unstructured MP without finetuning. 
When using Algorithm~\ref{alg:fm_callback}, 
we solved $\text{FM}(\mathcal{S})$ by setting
parameters MIPFocus, PoolSearchMode to 1 and PoolSolutions to 1000, ensuring that a greater number of solutions is obtained. 
We report the baseline of directly solving $\text{FM}(\mathcal{S})$ with default parameter values, 
since it gave better results in this case. 
We set the time limit to 600 seconds for all models. 
All other settings are as in Section~\ref{sec:nv_exp}.

\subsection{Results and Analysis}

Our approach based on Algorithm~\ref{alg:fm_callback} generally produced better results for the instances that would be typically regarded as harder to solve for having a larger number of linear regions~\cite{tong2024walk}. 
These are modest results in comparison with Section~\ref{sec:nv_exp}, but  
with noteworthy similarities and differences.  
Table~\ref{tab:fm_dimensions} summarizes the number of networks for which we found better values with Algorithm~\ref{alg:fm_callback}~by input size, number of layers (network depth), and neurons per layer (layer width).

We draw the following conclusions and justifications based on the results:
\begin{enumerate}[(A)]
\item The results from Algorithm~\ref{alg:fm_callback} are generally better if at least one of the dimensions of the neural network has a larger value (a),  
and consistently more so in the case of layer width (b).
\begin{enumerate}[(a)]
    \item In 32 out of 35 cases where any  dimension is not the smallest, at least half of the instances have a better solution. In contrast, larger input sizes were less favorable to our approach in the case of network verification.
    \item Unlike input size and network depth, any larger value for network width corresponds to a favorable case for Algorithm~\ref{alg:fm_callback}. In contrast, network width was the least relevant dimension in network verification. 
\end{enumerate}

\item Increasing along a dimension or along the pruning rate while fixing the other does not necessarily lead to better results (c), but very large networks yield considerably better results when using the highest pruning rate (d).
\begin{enumerate}[(a)]
    \setcounter{enumii}{2}
    \item We see at least non-monotonic variation of percentages if we increase any dimension along the same pruning rate, or if we increase the pruning rate along the same dimension. 
    In contrast, we generally observe a monotonic improvement up to 90\% pruning rate in network verification. 
    \item In all cases in which one of the dimensions is the largest, the best results are obtained for the highest pruning rate, 
    unlike in network verification.
\end{enumerate}
\end{enumerate}

\setlength{\tabcolsep}{9pt}
\begin{table}[h!]
\centering
\caption{Percentage of instances for which solving the function maximization problem indirectly on a pruned network yields a better solution by   
pruning rate (columns); 
disaggregated in terms of input size (top rows), number of layers (middle), and neurons per layer (bottom rows). 
Favorable figures (above 50\%) are reported in bold.
}\label{tab:fm_dimensions}
\begin{tabular}{l c || c c c c c}
  && \multicolumn{5}{c}{\textbf{Pruning Rate}} \\
  && 0.3 & 0.5 & 0.8 & 0.9 & 0.95 \\ 
  \hline \hline
  \multirow{3}{*}{\textbf{Input Size}} & 100 & 33.3\% & 40.0\% & 48.3\% & 50.0\% & 45.0\% \\
   & 1,000 & \textbf{63.3\%} & 48.3\% & \textbf{63.3\%} & 50.0\% & 50.0\% \\
   & 10,000 & \textbf{65.0\%} & \textbf{58.3\%} & \textbf{56.7\%} & \textbf{55.0\%} & \textbf{78.3\%} \\
   \hline
   \multirow{4}{*}{\textbf{Network Depth}} & 2 & 48.9\% & 40.0\% & 24.4\% & 24.4\% & 24.4\% \\
   & 3 & \textbf{57.8\%} & 44.4\% & \textbf{66.7\%} & \textbf{55.6\%} & \textbf{51.1\%} \\
   & 4 & \textbf{62.2\%} & \textbf{55.6\%} & \textbf{66.7\%} & \textbf{64.4\%} & \textbf{73.3\%} \\
   & 5 & 46.7\% & \textbf{55.6\%} & \textbf{66.7\%} & \textbf{62.2\%} & \textbf{82.2\%} \\
   \hline
   \multirow{3}{*}{\textbf{Layer Width}} & 50 & \textbf{51.7\%} & 43.3\% & 35.0\% & 31.7\% & 33.3\% \\
   & 100 & \textbf{58.3\%} & \textbf{51.7\%} & \textbf{65.0\%} & \textbf{55.0\%} & \textbf{60.0\%} \\
   & 200 & \textbf{51.7\%} & \textbf{51.7\%} & \textbf{68.3\%} & \textbf{68.3\%} & \textbf{80.0\%}
\end{tabular}
\end{table}


\section{Conclusion}

With the goal of solving optimization problems embedding a dense neural network, 
we tackled those problems indirectly through  drastically sparsified neural networks serving as surrogates. 
Those problems become very difficult to solve as the neural networks grow larger in size, 
and the surrogate is a pruned version of the same neural network. 
By making the models sparser, we naturally expect to find solutions faster, 
and in some cases we do not expect or do not need to necessarily find an optimal solution. 
We believe that this work contributes to understanding how to tackle constraint learning models more effectively.

We have found that a cost-effective approach is 
applying unstructured pruning while carefully choosing which connections to prune  
but not finetuning the pruned network afterwards.
We obtained consistently strong results in network verification, 
even if the surrogate network had very low accuracy. 
We also obtained encouraging, albeit modest, results in function maximization, 
but under different conditions than we found them for network verification. 

%

\begin{credits}
\subsubsection{\ackname} The authors were supported by the National Science Foundation grant 2104583, including Jiatai Tong while at Bucknell University.

\subsubsection{\discintname}
The authors have no competing interests to declare that are
relevant to the content of this article. 
\end{credits}
%
%
%
\bibliographystyle{splncs04}
\bibliography{references}

\end{document}